%% file: Hypoelliptic-damped-waves.tex
\documentclass[10pt]{article}

\usepackage[english,activeacute]{babel}
\usepackage{epsfig}

\usepackage[latin1]{inputenc}
\usepackage[T1]{fontenc}
\usepackage{stmaryrd}
\usepackage{amsmath,amsfonts,amssymb,mathrsfs,amsthm}
\usepackage{xcolor}
\usepackage{dsfont}
\usepackage{graphicx}
\usepackage[mathscr]{eucal}
\usepackage{makeidx}
\usepackage{verbatim}
\usepackage{graphics,graphicx}
\usepackage{textcomp}
\usepackage{float}
\usepackage[colorlinks=true]{hyperref}

\input mymacros.tex

\newcommand\bna{\begin{eqnarray*}} 
\newcommand\ena{\end{eqnarray*}}

\newcommand\bnan{\begin{eqnarray}} 
\newcommand\enan{\end{eqnarray}}

\newcommand\bnp{\begin{proof}} 
\newcommand\enp{\end{proof}}

\newcommand\bneq{\begin{eqnarray*}\left\lbrace \begin{array}{rcl}}
\newcommand\eneq{\end{array} \right.\end{eqnarray*}}
\newcommand\bneqn{\begin{eqnarray}\left\lbrace \begin{array}{rcl}}
\newcommand\eneqn{\end{array} \right.\end{eqnarray}}

\newcommand\nor[2]{\left\|#1\right\|_{#2}}

\renewcommand{\kappa}{\nu}

\numberwithin{equation}{section}

\newtheorem{hypo}{Assumption}[section]

\author{Camille Laurent\footnote{CNRS UMR 7598 and Sorbonne Universit\'es UPMC Univ Paris 06, Laboratoire Jacques-Louis Lions, F-75005, Paris, France, email: laurent@ann.jussieu.fr} and Matthieu L\'eautaud\footnote{Laboratoire de Math\'ematiques d'Orsay, Universit\'e Paris-Sud, CNRS, Universit\'e Paris-Saclay, B\^atiment 307, 91405 Orsay Cedex France, email: matthieu.leautaud@math.u-psud.fr.}
}

\begin{document}
\title{Logarithmic decay for damped hypoelliptic wave and Schr\"odinger equations}
\maketitle

\begin{abstract}
We consider damped wave (resp. Schr\"odinger and plate) equations driven by a hypoelliptic ``sum of squares'' operator $\scrL$ on a compact manifold and a damping function $b(x)$. We assume the Chow-Rashevski-H\"ormander condition at rank $k$ (at most $k$ Lie brackets needed to span the tangent space) together with analyticity of $\M$ and the coefficients of $\scrL$. We prove decay of the energy at rate $\log(t)^{-\frac{1}{k}}$ (resp. $\log(t)^{-\frac{2}{k}}$ ) for data in the domain of the generator of the associated group. We show that this decay is optimal on a family of Grushin-type operators. 
This result follows from a perturbative argument (of independent interest) showing, in a general abstract setting, that quantitative approximate observability/controllability results for wave-type equations imply {\em a priori} decay rates for associated damped wave, Schr\"odinger and plate equations. 
The adapted quantitative approximate observability/controllability theorem for hypoelliptic waves is obtained by the authors in~\cite{LL:15,LL:17Hypo}. 
\end{abstract}


%
%

\section{Introduction and statements}
\subsection{Damped hypoelliptic evolution equations}
\label{s:intro-hypo}
We consider a smooth compact connected $d$-dimensional manifold $\M$, endowed with a smooth positive density measure $ds$. We denote by $L^2 =L^2(\M) = L^2(\M, ds)$ the space of square integrable functions with respect to this measure. Given a smooth vector field $X$, we define by $X^*$ its formal dual operator for the duality of $L^2(\M)$, that is, 
 $$
 \int_\M X^*(u)(x) \overline{v(x)}ds(x) =   \int_\M u(x) \overline{X(v)(x)}ds(x), \quad \text{ for any } u, v \in C^\infty(\M) .
 $$
Given $m\in \N$ and $m$ smooth real vector fields $X_1,\cdots,X_m$, we consider the (H\"ormander's type~I) hypoelliptic operator (also called sub-Riemannian Laplacian, see e.g.~\cite[Remark~1.30]{LL:17Hypo})
\bnan
\label{def:L}
\scrL =   \sum_{i=1}^m X_i^*X_i .
\enan
Note that $\scrL$ is  formally symmetric nonnegative since $( \scrL u,v)_{L^2(\M)}=  \sum_{i=1}^m (X_i u , X_i v )_{L^2(\M)}$ for all $u,v \in C^\infty(\M)$.
Given a nonnegative (so-called damping) function $b \in L^\infty(\M ; \R_+)$, we are interested in the first place in asymptotic properties of the damped wave equation associated to $(\scrL,b)$
\begin{align}
\label{e:damped-hypo}
\begin{cases}
(\d_t^2 + \scrL  + b \d_t )u = 0, \quad \text{ on }(0, + \infty) \times \M ,\\
(u, \d_t u)|_{t=0} = (u_0 , u_1) , \quad \text{ on }\M .
\end{cases}
\end{align}
Solutions of~\eqref{e:damped-hypo} enjoy formally the following dissipation identity (obtained by taking the inner product of~\eqref{e:damped-hypo} with $\d_t u$ and integrating on $(0,T)$):
$$
E(u(T)) - E(u(0)) = -\int_0^T \int_\M b(x)|\d_t u (t,x)|^2 ds(x) \, dt ,\quad E(u) = \frac12 \left(  \sum_{i=1}^m \nor{X_i u}{L^2(\M)}^2+  \nor{\d_tu}{L^2(\M)}^2\right) .
$$
We are also interested in the damped Schr\"odinger equation associated to $(\scrL,b)$
\begin{align}
\label{e:damped-schro-hypo}
\begin{cases}
( i \d_t + \scrL  + i b )u = 0, \quad \text{ on }(0, + \infty) \times \M ,\\
u|_{t=0} = u_0 , \quad \text{ on }\M , 
\end{cases}
\end{align}
for which the $L^2$ norm is a dissipated quantity (obtained by taking imaginary part of the inner product of~\eqref{e:damped-schro-hypo} with $u$ and integrating on $(0,T)$):
$$
\frac12 \nor{u(T)}{L^2(\M)}^2 - \frac12\nor{u_0}{L^2(\M)}^2 = -\int_0^T \int_\M b(x)|u(t,x)|^2ds(x) \, dt .
$$
Another related equation with similar behavior is the damped plate equation associated to $(\scrL,b)$
\begin{align}
\label{e:damped-beam-hypo}
\begin{cases}
(\d_t^2 + \scrL^2  + b \d_t )u = 0, \quad \text{ on }(0, + \infty) \times \M ,\\
(u, \d_t u)|_{t=0} = (u_0 , u_1) , \quad \text{ on }\M .
\end{cases}
\end{align}
Solutions of~\eqref{e:damped-beam-hypo} also enjoy formally a similar dissipation identity
$$
E_P(u(T)) - E_P(u(0)) = -\int_0^T \int_\M b(x)|\d_t u (t,x)|^2 ds(x) \, dt ,\quad E_P(u) = \frac12 \left( \nor{\scrL u}{L^2(\M)}^2+  \nor{\d_tu}{L^2(\M)}^2\right) .
$$

Hence, in the three situations, ``energy'' decays, and an interesting question is to understand if it converges to zero, and if so, at which rate.

\medskip
We shall always assume throughout that the family $(X_i)$ satisfies the Chow-Rashevski-H\"ormander condition (or is ``bracket generating'').
\begin{hypo}
\label{assumLiek}
There exists $\ell \geq 1$ so that for any $x\in \M$, $\Lie^\ell(X_1,\cdots,X_m) (x)=T_x\M$. Denote then by $k \in \N^*$ the minimal $\ell$ for which this holds.
\end{hypo}
Here, $\Lie^\ell$ denotes the Lie algebra at rank $\ell$ of the vector fields. The integer $k$ is sometimes referred to as the \textit{hypoellipticity index} of $\scrL$.
Under Assumption~\ref{assumLiek}, the celebrated H\"ormander \cite{Hor67} and Rothschild-Stein \cite{RS:76} theorems (see also~\cite{BCN:82} for a simpler proof) state that $\scrL$ is 
subelliptic of order $\frac{1}{k}$, that is: there is $C>0$ such that for any $u\in C^{\infty}(\M)$, we have
\bnan
\nor{u}{H^{\frac{2}{k}}(\M)}^2\leq C\nor{\scrL u}{L^2(\M)}^2+C\nor{u}{L^2(\M)}^2 . \label{estimhypo3}
\enan
As a consequence, the operator $\scrL$ is selfadjoint on $L^2(\M)$ with domain $\scrL : D(\scrL) \subset L^2(\M) \to L^2(\M)$. Since $H^2(\M) \subset D(\scrL) \subset H^{\frac{2}{k}}(\M)$, $\scrL$ has compact resolvent and thus admits a Hilbert basis of eigenfunctions $(\varphi_j)_{j \in \N}$, associated with the real eigenvalues $(\lambda_j)_{j \in \N}$, sorted increasingly, that is
\bnan
\label{e:spectral-elts}
\scrL\varphi_i=\lambda_i\varphi_i, \qquad (\varphi_i , \varphi_j)_{L^2(\M)} = \delta_{ij} , \qquad 0= \lambda_0 < \lambda_1 \leq \lambda_2 \leq \cdots \leq \lambda_j \to + \infty .
\enan
This allows in particular to define adapted Sobolev spaces:  
$$
 \H^s_\scrL =  \{ u \in \D'(\M) ,  \ (1+\scrL)^\frac{s}{2}   u \in L^2(\M) \} , \quad \nor{u}{\H^s_\scrL}=\nor{(1+\scrL)^\frac{s}{2} u}{L^2(\M)} , \quad  s\in \R ,
$$
where $f(\scrL) u = \sum_{j \in \N} f(\lambda_j) (u , \varphi_j)_{L^2(\M)} \varphi_j$.

In addition to Assumption~\ref{assumLiek}, we will also make the following analyticity assumption. 
\begin{hypo}
\label{hypoanal}
The manifold $\M$, the density $ds$, and the vector fields $X_i$ are real-analytic. 
\end{hypo}
A non-exhaustive list of classical examples of operators $\scrL$ encompassed by this frameworks is provided in~\cite[Section~1.1]{LL:17Hypo}.
Note that the damping function $b$ does not need to be analytic but only $L^\infty$; in particular our results work for $b=\mathds{1}_{\omega}$ if $\omega$  is a non-empty open subset of $\M$.

On the space $ \H^1_\scrL \times L^2$, the operator $\A = 
\left(
\begin{array}{cc}
0   &  \id \\
- \scrL & - b(x)
\end{array}
\right)$ with $D(\A) = \H^2_\scrL \times \H^1_\scrL$ generates a bounded semigroup (from the Hille-Yosida theorem) and~\eqref{e:damped-hypo} admits a unique solution $u \in C^0(\R^+ ; \H^1_\scrL )\cap C^1(\R^+ ;L^2)$. Our main results for damped hypoelliptic waves are summarized in the following two theorems.
\begin{theorem}[Decay rates for damped hypoelliptic waves]
\label{e:decay}
Assume that $b\in L^\infty(\M)$ is such that $b\geq \delta>0$ a.e. on a nonempty open set, together with Assumptions~\ref{assumLiek} and~\ref{hypoanal}. Then, for all $(u_0, u_1) \in \H^1_\scrL \times L^2$, the associated solution to~\eqref{e:damped-hypo} satisfies $E(u(t)) \to 0$. Moreover, for all $j \in \N^*$, there exists $C_j>0$ such that for all $(u_0, u_1) \in D(\A^j)$, the associated solution to~\eqref{e:damped-hypo} satisfies 
\begin{equation}
\label{e:log-decay}
E(u(t))^{\frac12} \leq \frac{C_j}{\log(t+2)^{j/k}} \nor{\A^j (u_0, u_1)}{\H^1_\scrL \times L^2} , \quad \text{ for all } t\geq 0 .
\end{equation}
\end{theorem}

Theorem~\ref{e:decay} is actually a consequence of the following result, together with~\cite{BD:08}. 
\begin{theorem}[Spectral properties for damped hypoelliptic waves]
\label{t:spectral}
Assume that $b\geq \delta>0$ a.e. on a nonempty open set, together with Assumptions~\ref{assumLiek} and~\ref{hypoanal}. Then, the spectrum of $\A$ contains only isolated eigenvalues with finite multiplicity, and satisfies:
\begin{enumerate}
\item \label{truc} $\ovl{\Sp(\A)}=\Sp(\A)$ and $\ker(\A) = \vect\{(1,0)\}$ (where $1$ denotes the constant function),
\item \label{i:loc-spec} $
\Sp(\A) \subset \left( \big( - \frac12 \|b\|_{L^\infty(\M)} , 0 \big) + i\R\right) \cup 
\left([ -  \|b\|_{L^\infty(\M)} , 0] + 0 i \right),$
\item \label{e:expo-bound} there exist $C,\kappa >0$ such that  $\nor{(is-\A)^{-1}}{\L(\H^1_\scrL \times L^2)} \leq C e^{\kappa |s|^k}$ for all $|s|\geq1$,
\item there exist $\eps ,\kappa >0$ such that $\Sp(\A) \cap \Gamma_k(\eps,\kappa ) = \{0\}$, where $\Gamma_k(\eps,\kappa ) = \{z \in \C , \Re(z) \geq - \eps e^{-\kappa |\Im(z)|^k}\}$.
\end{enumerate}
\end{theorem}
The first two points are rather standard, see~\cite{Leb:96}. Item~\ref{e:expo-bound} is the key information in the Theorem, and is a consequence of the main theorem in~\cite[Theorem~1.15]{LL:17Hypo}. 
The last point of the theorem states an exponentially small spectral gap, and is a consequence of Item~\ref{e:expo-bound} and Neumann series expansion.

Combined together, Theorems~\ref{e:decay} and~\ref{t:spectral} are the counterparts to~\cite[Th\'eor\`eme~1]{Leb:96} in the case of the usual wave equation ($k=1$, in which case no analyticity is required, and boundary condition can be dealt with).

Note that the fact that $\Sp(\A)\cap i\R = \{0\}$ in Item~\ref{i:loc-spec} (which, in turn, implies that $E(u(t)) \to 0$ in Theorem~\eqref{e:decay} for all solutions to~\eqref{e:damped-hypo}) is actually a consequence of the {\em qualitative} uniqueness: 
\bnan
\label{e:UCP-eig-Hypo}
\Big( \varphi \in \H_\scrL^2 , \quad z \in \C, \quad \scrL \varphi  = z \varphi \text{ on }\M , \quad \varphi = 0 \text{ on } \omega \Big) \implies \varphi \equiv 0\text{ on }\M ,
\enan
 proved by Bony~\cite{Bo:69}, as a consequence of the Holmgren-John theorem. Even this weaker property is not well understood for general hypoelliptic operators if we drop Assumption~\ref{hypoanal}, see~\cite{Ba:86}.
Here the key point is the quantification of the Holmgren-John theorem proved in~\cite{LL:15,LL:17Hypo}.

\bigskip
We present analogue results in the case of the damped hypoelliptic Schr\"odinger equation.
We set $\A_S := i \scrL - b$ with $D(\A_S)=D(\scrL)$, so that~\eqref{e:damped-schro-hypo} reformulates as $(\d_t- \A_S )u=0$. Note that $\A_S$ generates a contraction semigroup (from the Hille-Yosida theorem) and~\eqref{e:damped-schro-hypo} admits a unique solution $u \in C^0(\R^+ ; L^2(\M) )$. 
Our main results for the damped hypoelliptic Schr\"odinger equation are summarized in the following two theorems.
\begin{theorem}[Decay rates for the damped hypoelliptic Schr\"odinger equation]
\label{e:decay-schro}
Assume that $b\in L^\infty(\M)$ is such that $b\geq \delta>0$ a.e. on a nonempty open set, together with Assumptions~\ref{assumLiek} and~\ref{hypoanal}. Then, for all $u_0 \in L^2(\M)$, the associated solution to~\eqref{e:damped-schro-hypo} satisfies $u(t) \to 0$ in $L^2(\M)$. Moreover, for all $j \in \N^*$, there exists $C_j>0$ such that for all $u_0\in D(\A_S^j)$, the associated solution to~\eqref{e:damped-schro-hypo} satisfies 
\begin{equation}
\label{e:log-decay-schro}
\nor{u(t)}{L^2(\M)} \leq \frac{C_j}{\log(t+2)^{2j/k}} \nor{\A_S^j u}{L^2(\M)} , \quad \text{ for all } t\geq 0 .
\end{equation}
\end{theorem}
Note that when comparing~\eqref{e:log-decay-schro} to~\eqref{e:log-decay}, the decay looks better ($\log(t+2)^{-2j/k}$ instead of $\log(t+2)^{-j/k}$) but actually consumes more derivatives: for smooth $b$, $\nor{\A_S^j u}{L^2(\M)} \simeq \nor{u}{\H_\scrL^{2j}}$ whereas $\nor{\A^j U}{L^2(\M)} \simeq \nor{U_0}{\H_\scrL^{j} \times \H_\scrL^{j-1}}$. Hence both decay rates essentially coincide for data having the same regularity.
Theorem~\ref{e:decay-schro} is a consequence of the following result, together with~\cite{BD:08}. 
\begin{theorem}[Spectral properties for the damped hypoelliptic Schr\"odinger equation]
\label{t:spectral-schro}
Assume that $b\geq \delta>0$ a.e. on a nonempty open set, together with Assumptions~\ref{assumLiek} and~\ref{hypoanal}. Then, the spectrum of $\A_S$ contains only isolated eigenvalues with finite multiplicity, and satisfies:
\begin{enumerate}
\item \label{i:loc-spec-schro} $
\Sp(\A_S) \subset \big[ -  \|b\|_{L^\infty(\M)} , 0 \big) + i[0,+\infty)$,
\item \label{e:expo-bound-schro} there exist $C,\kappa >0$ such that  $\nor{(is-\A_S)^{-1}}{\L(L^2)} \leq C e^{\kappa |s|^{k/2}}$ for all $s\in \R$,
\item there exist $\eps ,\kappa >0$ such that $\Sp(\A_S) \cap \Gamma_{k,S}(\eps,\kappa ) =\emptyset$, where $\Gamma_{k,S}(\eps,\kappa ) = \{z \in \C , \Re(z) \geq - \eps e^{-\kappa |\Im(z)|^{k/2}}\}$.
\end{enumerate}
\end{theorem}
Note that in the elliptic case $k=1$, the results of Theorems~\ref{e:decay-schro},~\ref{t:spectral-schro} are more or less classical, even though we did not see them written explicitely in the literature. In this situation, analyticity is not necessary and boundary value problems can be dealt with. Our abstract perturbative proof below works as well, as a consequence of~\cite{LL:15} (with Dirichlet boundary conditions). One can however start from the seminal estimates of Lebeau-Robbiano in this situation, see~\cite{LR:95,Leb:96} for Dirichlet conditions and~\cite{LR:97} for Neumann boundary conditions.

\bigskip
A similar result holds for the plate equation. The framework is quite similar to the wave equation. We will work on the space $ \H^2_\scrL \times L^2$ with the operator $\A_P = 
\left(
\begin{array}{cc}
0   &  \id \\
- \scrL^2 & - b(x)
\end{array}
\right)$ with $D(\A_P) = \H^4_\scrL \times \H^2_\scrL$. It generates a bounded semigroup and~\eqref{e:damped-beam-hypo} admits a unique solution $u \in C^0(\R^+ ; \H^2_\scrL )\cap C^1(\R^+ ;L^2)$.
\begin{theorem}[Decay rates for damped hypoelliptic plates]
\label{e:decaybeam}
Assume that $b\in L^\infty(\M)$ is such that $b\geq \delta>0$ a.e. on a nonempty open set, together with Assumptions~\ref{assumLiek} and~\ref{hypoanal}. Then, for all $(u_0, u_1) \in \H^2_\scrL \times L^2$, the associated solution to~\eqref{e:damped-beam-hypo} satisfies $E_P(u(t)) \to 0$. Moreover, for all $j \in \N^*$, there exists $C_j>0$ such that for all $(u_0, u_1) \in D(\A^j)$, the associated solution to~\eqref{e:damped-beam-hypo} satisfies 
\begin{equation}
\label{e:log-decaybeam}
E_P(u(t))^{\frac12} \leq \frac{C_j}{\log(t+2)^{2j/k}} \nor{\A_P^j (u_0, u_1)}{\H^2_\scrL \times L^2} , \quad \text{ for all } t\geq 0 .
\end{equation}
\end{theorem}
Similar spectral statements as Theorems~\ref{t:spectral} and~\ref{t:spectral-schro} hold for the plate equation. We leave the details to the reader. Again, by using our result~\cite{LL:15}, we could also obtain a logarithmic decay in the elliptic case $k=1$ for a compact manifold with boundary and with Dirichlet boundary conditions. We do not know if this result is new in this context. The literature is quite big, we refer to \cite{Leb:92} and \cite{K:92} for exact control results (implying exponential decay of the damped equation) and e.g. to~\cite{ADZ:14} for a spectral analysis of the decay rate.

\bigskip
Finally, we show that results of Theorems~\ref{e:decay},~\ref{t:spectral},~\ref{e:decay-schro},~\ref{t:spectral-schro} are optimal in general (in case $k>1$; this is already known in the elliptic case $k=1$, see~\cite{Leb:96,LR:97}). This is also the case for Theorem~\ref{e:decaybeam} (and the associated spectral statement); we do not state the result for the sake of brevity.
\begin{proposition}
\label{Prop:BCG-damped}
Consider the manifold with boundary $\M=[-1,1]\times(\R/\Z)$, endowed with the Lebesgue measure $dx$, and for $k \in (1, +\infty)$, define the operator $\scrL =  - \big( \partial_{x_1}^2 + x_1^{2(k-1)}\partial_{x_2}^2 \big)$, with Dirichlet conditions on $\d \M$.  
Assume that $\supp(b) \cap \left\{x_1=0\right\} = \emptyset$. Then there exist  $C, \kappa>0$ and a sequence $(s_j)_{j\in \N}$ with $s_j\to + \infty$  such that 
\begin{align}
\label{e:lower-resolvent-bounds}
\nor{(is_j-\A)^{-1}}{\L(\H^1_\scrL \times L^2)} \geq C e^{\kappa s_j^k} , \quad  \nor{(is_j-\A_S)^{-1}}{\L(\H^1_\scrL \times L^2)} \geq C e^{\kappa s_j^{k/2}} , \quad  \text{for all } j \in \N.
\end{align}
Moreover, if for all $(u_0, u_1) \in D(\A)$, the associated solution to~\eqref{e:damped-hypo} satisfies
\begin{equation*}
E(u(t))^{\frac12} \leq  f(t) \nor{\A(u_0, u_1)}{\H^1_\scrL \times L^2} , \quad \text{ for all }  t\geq 2 ,
\end{equation*}
then there is $C>0$ such that $f(t) \geq \frac{C}{\log(t)^{1/k}}$. Similarly, if for all $u_0 \in \H_\scrL^1$, the associated solution to~\eqref{e:damped-schro-hypo} satisfies
\begin{equation*}
\nor{u(t)}{L^2(\M)} \leq f(t)\nor{\A_S u}{L^2(\M)} , \quad \text{ for all } t\geq 2 ,
\end{equation*}
then there is $C>0$ such that $f(t) \geq \frac{C}{\log(t)^{2/k}}$.
\end{proposition}
Recall that for $k\in \N^*$, the operator $\scrL =  - \big( \partial_{x_1}^2 + x_1^{2(k-1)}\partial_{x_2}^2 \big)$ satisfies precisely Assumption~\ref{assumLiek}.
The first part is thus a consequence of \cite[Section~2.3]{BeauchardCanGugl} as reformulated in~\cite[Proposition~1.14]{LL:17Hypo}. It prove the optimality in general of Item~\ref{t:spectral} in Theorem~\ref{e:expo-bound}.
The second part is a corollary of the first together with~\cite{BD:08}, and proves optimality of~\eqref{e:log-decay}.

A reformulation of~\cite{Letrouit:20} (e.g. together with~\cite{Har:89}) in the present context states that if $\vect (X_1(x), \cdots , X_m(x)) \neq T_x\M$ for $x$ in a dense subset of $\M$, and $\M \setminus \supp(b) \neq \emptyset$, then uniform decay does not hold: there is no function $f : \R^+ \to \R^+$ with $f(t)\to 0$ such that $E(u(t)) \leq f(t)E(u(0))$. This contrasts with the Riemannian case~\cite{RT:74,BLR:92}, and gives in this context a stronger interest to the result of Theorem~\ref{e:decay} as compared to the Riemannian counterpart.

However, one may notice that logarithmic decay as in Theorem~\ref{e:decay} is not always optimal. Combining for instance~\cite[Theorem~1]{BS:19} together with~\cite[Theorem~2.3]{AL:14} implies that $\frac{C_j}{\log(t+2)^{j/k}}$ in~\eqref{e:log-decay} can be replaced by $\frac{C_j}{t^{j/2}}$ (and this is probably not optimal) in the geometric setting of Proposition~\ref{Prop:BCG-damped} if $b(x_1,x_2) = \mathds{1}_{(a,b)}(x_2)$, for any $a<b$.

Similarly, logarithmic decay in Theorem~\ref{e:decay-schro} is not always optimal. For instance~\cite[Theorem~1]{BS:19} (together with classical equivalence between observability for the conservative system and uniform stabilization for the damped system) implies that in the geometric setting of Proposition~\ref{Prop:BCG-damped} if $b(x_1,x_2) = \mathds{1}_{(a,b)}(x_2)$ for $a<b$, then uniform decay holds, that is: there are $C,\gamma>0$ such that $\nor{u(t)}{L^2} \leq C e^{-\gamma t}\nor{u_0}{L^2}$ for all solutions to~\eqref{e:damped-schro-hypo}.

\bigskip
Let us finally remark that all proofs below rely on the approximate observability/controllability of the hypoelliptic wave equation with optimal cost. The latter result is proved by the authors in~\cite{LL:17Hypo}.
It is interesting to notice that in the elliptic case ($k=1$ in the discussion above), the approximate observability/controllability of the wave equation (proved in~\cite{LL:15}) with optimal (exponential) cost allows to recover many known control results obtained with Carleman estimates.
In particular, it implies 
\begin{enumerate}
\item null-controllability of the heat equation with optimal short-time behavior, as proved in~\cite{EZ:11} and~\cite[Proposition~1.7]{LL:18} (the original result is~\cite{LR:95,FI:96}),
\item approximate observability/controllability of the heat equation with optimal (exponential) cost~\cite[Chapter~4]{LL:17Hypo} (the original result is~\cite{FCZ:00}),
\item optimal logarithmic decay for the damped wave equation, see Theorem \ref{e:decay} for $k=1$ (the original result is~\cite{Leb:96,LR:97}).
\end{enumerate}
Here, we provide a proof of the last point in a general framework presented in Section~\ref{s:abstract} below, and deduce counterparts for hypoelliptic equations using~\cite{LL:17Hypo}.

\subsection{From approximate control to damped waves : abstract setting}
\label{s:abstract}

As already mentioned, we prove the results of Theorems~\ref{t:spectral} and~\ref{e:decay} in an abstract operator setting. This allows us to stress links between the cost of approximate controls and a priori decay rates for damped waves.
This follows the spirit of e.g.~\cite{Har:89,BZ:04,Phung:01,Miller:05,Miller:06b,TW:09,EZ:11,AL:14,CPSST:19},
exploring the links between different equations and their control
properties (i.e. observability, controllability,
stabilization...). 
Here, we follow closely~\cite{AL:14}. 

Let $H$ and $Y$ be two Hilbert spaces (resp. the state space and the observation/control space) with norms $\| \cdot \|_H$ and $\| \cdot \|_Y$, and associated inner products $( \cdot , \cdot )_H$ and $( \cdot , \cdot )_Y$. 
We denote by $A : D(A)\subset H \to H$ a {\em nonnegative} selfadjoint operator with compact resolvent, and $B \in \L(Y;H)$ a control operator. We recall that $B^* \in \L(H;Y)$ is defined by $(B^* h ,y)_Y = (h, B y)_H$ for all $h \in H$ and $y \in Y$. We define $H_1 = D(A^\frac12)$, equipped with the graph norm $\nor{u}{H_1} := \|(A+\id)^\frac12 u \|_H$, and its dual $H_{-1} = (H_1)'$ (using $H$ as a pivot space) endowed with the norm $\nor{u}{H_{-1}}:=\|(A+\id)^{-\frac12} u \|_H$.

In applications to Theorems~\ref{e:decay}-\ref{t:spectral}-\ref{e:decay-schro}-\ref{t:spectral-schro}, we take $H=Y=L^2(\M)$, $A=\scrL$ and $B = B^*$ is multiplication by the function $\sqrt{b}$.

We introduce in this abstract setting the wave equation
\begin{align}
\label{eq: abstract waves}
\begin{cases}
\d_t^2 u + A u =  F  , \\
(u, \d_t u)|_{t=0} = (u_0 , u_1) , 
\end{cases}
\end{align}
the damped wave equation
\begin{align}
\label{eq: damped abstract waves}
\begin{cases}
\d_t^2 u + A u + B B^* \d_t u = 0, \\
(u, \d_t u)|_{t=0} = (u_0 , u_1) , 
\end{cases}
\end{align}
and the damped Schr\"odinger equation
\begin{align}
\label{eq: damped abstract schro}
\begin{cases}
i \d_t u + A u + i B B^* u = 0, \\
u|_{t=0} = u_0 .
\end{cases}
\end{align}
 
\begin{definition}
\label{d:approx-obs}
Given $T>0$ and a function $G : \R_+ \to \R_+$, we say that the wave equation~\eqref{eq: abstract waves} with $F=0$ is approximately observable from $B^*$ in time $T$ with cost $G$ if there is $\mu_0>0$ such that for all $(u_0,u_1)\in H_1 \times H$, the associated solution $u$ to~\eqref{eq: abstract waves} with $F=0$ satisfies 
\bnan
\label{th-estimate-k-bis}
\nor{(u_0,u_1)}{H\times H_{-1}} \leq G(\mu) \nor{B^*u}{L^2(0,T;Y)} +\frac{1}{\mu}\nor{(u_0,u_1)}{H_1 \times H} , \quad \text{for all }\mu \geq \mu_0 .
\enan
\end{definition}
According to~\cite{Robbiano:95} or~\cite[Appendix]{LL:17approx}, this is equivalent to approximate controllability ($\eps$ close) with cost $G(1/\eps)$.
This is satisfied for the usual wave equation in a general context with $B^* = \mathds{1}_\omega$, $G(\mu) = Ce^{\kappa \mu}$, for all $T> 2\sup_{x\in \M} d_g(x,\omega)$ (where $d_g$ is the Riemannian distance), as proved in~\cite{LL:15}.
For the hypoelliptic wave equation, we proved in \cite[Theorem~1.15]{LL:17Hypo} that this is satisfied for $B^* = \mathds{1}_\omega$, $G(\mu) = Ce^{\kappa \mu^k}$, for all $T> 2\sup_{x\in \M} d_\scrL(x,\omega)$ (where $d_\scrL$ is the appropriate sub-Riemannian distance and $k$ the hypoellipticity index of $\scrL$).

\bigskip
Our main results can be divided in several steps. Firstly we have
\begin{proposition} 
\label{c:res-fct-propres}
Let $G : \R_+ \to \R_+$ be such that $G(\mu) \geq \frac{c_0}{\mu}>0$ for  $\mu \geq \mu_0$. Assume that there is $T>0$ such that the wave equation~\eqref{eq: abstract waves} with $F=0$ is approximately observable from $B^*$ in time $T$ with cost $G$ in the sense of Definition~\ref{d:approx-obs}. 
Then, we have 
\bnan
\label{e::UCPeigenvalue}
 \big(\lambda \in \C,\  v\in D(A),\ Av = \lambda^2 v ,\ B^*v=0 \big)\implies v=0 ,
\enan
and there is $\lambda_0>0$ such that for all $\alpha >0$,  
\begin{align}
\label{e:res-fct-prop}
\nor{v}{H}\leq \frac{K}{\alpha}(\lambda + \sqrt{2}+\alpha)  G(\lambda + \sqrt{2}+\alpha) \big( \nor{B^* v}{Y} + C\nor{(A - \lambda^2)v }{H}  \big) , \quad \text{ for all }v \in D(A) , \lambda \geq \lambda_0 .
\end{align}
with $K = \sqrt{T} + c_0^{-1}$ and $C>0$ a constant depending only on $B$ and $T$.
\end{proposition}
Note that in this statement, $\sqrt{2}$ can be replaced by $1$ at the cost of a slightly longer proof, and $\lambda_0$ is the $\mu_0$ given in the definition of approximate observability. In most applications we have in mind, $G(\mu) \approx e^{\kappa \mu^k}$ and the estimate is better for smaller values of $\alpha$. In a situations in which one would have $G(\mu) \approx \mu^\gamma$, then a better choice of $\alpha$ would be $\alpha \approx \lambda$, so that~\eqref{e:res-fct-prop} remains a bound of order $G(\lambda)$.

\medskip
Secondly, we assume that for some function $\mathsf{G}$ and some $\lambda_0>0$ we have
\begin{align}
\label{e:res-fct-prop-bis}
\nor{v}{H}\leq \mathsf{G}(\lambda) \big( \nor{B^* v}{Y} + \nor{(A - \lambda^2)v }{H}  \big) , \quad \text{ for all }v \in D(A) , \lambda \geq \lambda_0 .
\end{align}
This is precisely~\eqref{e:res-fct-prop} with $\mathsf{G}(\lambda) = \frac{K(1+C)}{\alpha}(\lambda + \sqrt{2}+\alpha)  G(\lambda + \sqrt{2}+\alpha)$.
From this estimate, we deduce the sought spectral properties (resolvent estimates and localization of the spectrum linked to the function $\mathsf{G}$). See Section~\ref{s:damped-schro} for the damped Schr\"odinger equation and Section~\ref{s:damped-wave} for the damped wave equation.
A direct application of Proposition~\ref{c:res-fct-propres} gives in the context of hypoelliptic operators. 
\begin{corollary}
\label{corresolvL}
With the notations of Section~\ref{s:intro-hypo}, assume that $b\in L^\infty(\M)$ is such that $b\geq \delta>0$ a.e. on a nonempty open set, together with Assumptions~\ref{assumLiek} and~\ref{hypoanal}.
Then,~\eqref{e:UCP-eig-Hypo} is satisfied and there is $\nu>0$, $C>0$ and $\lambda_0>0$ such that,  
\begin{align*}
\nor{v}{L^2(\M)}\leq C e^{\kappa\lambda^k}\big( \nor{b v}{L^2(\M)} + \nor{(\scrL - \lambda^2)v }{L^2(\M)}  \big) , \quad \text{ for all }v \in \H^2_{\scrL} , \lambda \geq \lambda_0 .
\end{align*}
\end{corollary}
This corollary states a stronger version of the Eigenfunction tunneling estimates of~\cite[Theorem~1.12]{LL:17Hypo} (which is the same statement for solutions to $(\scrL - \lambda^2)v=0$).
Note that the constant $\kappa$ is (essentially) the same as in the cost of approximate controls in~\cite[Theorem~1.15]{LL:17Hypo}.

\medskip
Thirdly, we deduce from the spectral properties the sought decay estimates (respectively in Sections~\ref{s:damped-schro} and~\ref{s:damped-wave} for the damped Schr\"odinger and wave equations) using the Batty-Duyckaerts theorem, which we now recall.
\begin{theorem}[Batty and Duyckaerts~\cite{BD:08}]
\label{t:batty-duyckaerts}
Let $(e^{t\B})_{t\geq 0}$ be a bounded $\Con^0$-semigroup on a Banach space $\mathcal{X}$, generated by $\B$. 

Assume that $ \nor{e^{t\B}\B^{-1}}{\L(\dot{\H})} \leq f(t)$ with $f$ decreasing to $0$. Then $i \R \cap \Sp(\B) = \emptyset$ and there is $C>0$ such that
$$
 \nor{(i\lambda-\B)^{-1}}{\L(\mathcal{X})} \leq 1 + C f^{-1} \left(\frac{1}{2(|\lambda|+1)} \right).
$$

Conversely, suppose that $i \R \cap \Sp(\B) = \emptyset$ and  
\begin{align}
\label{e:res-bounded-M}
 \nor{(is-\B)^{-1}}{\L(\mathcal{X})}  \leq \mathsf{M}(|s|) , \quad s \in \R  ,
\end{align}
where $\mathsf{M} : \R_+\to \R_+^*$ is a non-decreasing function on $\R_+$. Then, setting 
\begin{equation}
\label{Mlog}
\mathsf{M}_{\log}(s) =  \mathsf{M}(s) \big( \log(1+ \mathsf{M}(s)) + \log(1+ s)\big),
\end{equation}
 there exists $c>0$ such that for all $j \in \N$,
$$
 \nor{e^{t\B}\B^{-j}}{\L(\mathcal{X})}   = \O \bigg(  \frac{1}{\mathsf{M}_{\log}^{-1}\left( \frac{t}{c j} \right)^j}\bigg) , \quad \text{as } t \to +\infty ,
$$
where $\mathsf{M}_{\log}^{-1} : \R^+\to \R^+$ denotes the inverse of the strictly increasing function $\mathsf{M}_{\log}$.
\end{theorem}
We refer to~\cite{Duy:15,CS:16} for alternative proofs of the result of~\cite{BD:08}.
Note that on a Hilbert space (which is the case here) $\mathsf{M}_{\log}$ in the result can be replaced by $\mathsf{M}$ if it is polynomial at infinity, according to~\cite[Theorem 2.4]{BT:10} (see also~\cite{CPSST:19} and the references therein for generalizations of~\cite{BT:10}).

\bigskip
\noindent
{\em Acknowledgements.} 
The first author is partially supported by the Agence Nationale de la Recherche under grant  SRGI ANR-15-CE40-0018.
Both authors are partially supported by the Agence Nationale de la Recherche under grant ISDEEC ANR-16-CE40-0013.

\section{Proofs}
\label{s:proofs}
\subsection{From approximate observability of waves to a free resolvent estimate with an observation term: Proof of Proposition~\ref{c:res-fct-propres}}

From approximate observability, we deduce the following (seemingly more general) result, concerning equation~\eqref{eq: abstract waves} with a general right hand-side $F$.

\begin{proposition} 
\label{p:extrait-hypotri}
Let $T>0$ and a function $G : \R_+ \to \R_+$. Assume that the wave equation~\eqref{eq: abstract waves} with $F=0$ is approximately observable from $B^*$ in time $T$ with cost $G$, in the sense of Definition~\ref{d:approx-obs}. Then, there are $\mu_0,C>0$ such that for all $F\in L^2(0,T; H)$ and $(u_0,u_1)\in H_1 \times H$, the associated solution $u$ to~\eqref{eq: abstract waves} satisfies 
\bnan
\label{th-estimate-k-bis-rhs}
\nor{(u_0,u_1)}{H\times H_{-1}} \leq G(\mu) \big(  \nor{B^*u}{L^2(0,T;Y)}  + C \nor{F}{L^2(0,T;H)}\big) +\frac{1}{\mu}\nor{(u_0,u_1)}{H_1 \times H} , \quad \text{for all }\mu \geq \mu_0 .
\enan
\end{proposition}
Note that the constant $\mu_0$ is actually the same as in Definition~\ref{d:approx-obs} and that $C$ depends only on $T$ and $\nor{B^*}{\L(Y;H)}$.
\bnp
According to the linearity of~\eqref{eq: abstract waves}, we decompose $u$ as $u=u^0+u^F$ where $u^0$ is the solution to~\eqref{eq: abstract waves} for $F=0$ and $u^F$ is the solution to~\eqref{eq: abstract waves} with $(u_0,u_1)=(0,0)$.

First, according to the assumption, Definition~\ref{d:approx-obs} applies to the function $u^0$, so that~\eqref{th-estimate-k-bis} reads:
\bnan
\label{th-estimate-k-bis-ter}
\nor{(u_0,u_1)}{H\times H_{-1}} \leq G(\mu) \nor{B^*u^0}{L^2(0,T;Y)} +\frac{1}{\mu}\nor{(u_0,u_1)}{H_1 \times H} , \quad \text{for all }\mu \geq \mu_0 .
\enan
Second, to estimate $u^F$, we perform classical energy inequalities for~\eqref{eq: abstract waves}. We rewrite~\eqref{eq: abstract waves} as 
$$
(\d_t^2 + A+\id ) u^F = u^F + F , \quad (u^F, \d_t u^F)|_{t=0} = (0,0) .
$$
Taking the inner product of this equation with $\d_tu^F$ (assuming at first that $F \in L^1_{\loc}(\R; H_1)$ and thus $u^F \in  C^0(\R;D(A))\cap C^1(\R ; H_1)\cap C^2(\R ; H)$) implies
$$
\frac12 \frac{d}{dt} \left( \nor{\d_t u^F}{H}^2  + \nor{u^F}{H_1}^2 \right) \leq \left( \nor{u^F}{H} + \nor{F}{H} \right) \nor{\d_t u^F}{H}.
$$
Writing $\tilde{E}(t) = \frac12\left( \nor{\d_t u^F}{H}^2  + \nor{u^F}{H_1}^2 \right)$, this is $\tilde{E}'(t) \leq 2\tilde{E}(t) + \nor{F}{H}^2$. The Gronwall lemma together with the vanishing initial data imply
 $$
  \sup_{t\in [0,T]} \nor{u^F(t)}{H}^2 \leq \sup_{t\in [0,T]} \tilde{E}(t)  \leq C_T \nor{F}{L^1(0,T;H)}^2 .
 $$
 As a consequence, boundedness of $B^*$ yields
 $$
 \nor{B^*u^F}{L^2(0,T;Y)} \leq \nor{B^*}{\L(Y;H)}\nor{u^F}{L^2(0,T;H)} \leq \nor{B^*}{\L(Y;H)}C_T\nor{F}{L^2(0,T; H)} .
 $$
Recalling that $u^0 = u - u^F$  and combining this estimate with~\eqref{th-estimate-k-bis-ter} yields for all $\mu \geq \mu_0$
\begin{align*}
\nor{(u_0,u_1)}{H\times H_{-1}}&  \leq G(\mu) \nor{B^*(u-u^F)}{L^2(0,T;Y)} +\frac{1}{\mu}\nor{(u_0,u_1)}{H_1 \times H} \\
&  \leq G(\mu) \left( \nor{B^*u}{L^2(0,T;Y)} + C_{B,T}\nor{F}{L^2(0,T; H)} \right) +\frac{1}{\mu}\nor{(u_0,u_1)}{H_1 \times H} ,
\end{align*}
which concludes the proof of the proposition.
\enp

From this result, we deduce a proof of Proposition~\ref{c:res-fct-propres} as a direct corollary.

\bnp[Proof of Proposition~\ref{c:res-fct-propres}]
For $v \in D(A)$ and $\lambda\in \C$, we may apply the result of Proposition~\ref{p:extrait-hypotri} to the function $u(t) = \cos(\lambda t)v$ which satisfies~\eqref{eq: abstract waves} with 
$$
u_0 = v , \quad u_1 =0 , \quad F(t) = \cos(\lambda t) (-\lambda^2 + A)v .
$$
We remark that the assumption of~\eqref{e::UCPeigenvalue} implies $F=0$ and $B^*u=0$, and hence~\eqref{th-estimate-k-bis-rhs} reads $\nor{v}{H} \leq \frac{1}{\mu}\nor{v}{H_1}$ for all $\mu \geq \mu_0$. Letting $\mu$ converges to $+\infty$ yields the conclusion of~\eqref{e::UCPeigenvalue}.

\medskip
Let us now prove~\eqref{e:res-fct-prop}.  Still for $u(t) = \cos(\lambda t)v$, we have 
$$
 \nor{B^*u}{L^2(0,T;Y)}^2 \leq T \nor{B^*v}{Y}^2,  \quad \nor{F}{L^2(0,T;H)}^2 \leq T \nor{(-\lambda^2 + A)v}{H}^2.
$$
Estimate~\eqref{th-estimate-k-bis-rhs} thus implies for all $\lambda \geq 0$, $\mu \geq \mu_0$
\begin{align}
\label{e:interm-interm}
\nor{v}{H}\leq G(\mu) \sqrt{T} \big( \nor{B^*v}{Y} + C\nor{(A - \lambda^2)v }{H}  \big) +\frac{1}{\mu}  \nor{v}{H_1}  .
\end{align}
We now remark that
$$
 \left( A v ,v \right)_{H} -  \lambda^2 \nor{v}{H}^2 = \left( (A - \lambda^2)v ,v \right)_{H}
 \leq \nor{(A - \lambda^2)v }{H}\nor{v}{H}
$$
Hence, we deduce
\begin{align*}
\nor{v}{H_1}^2 & =  \left( (A +1)v ,v \right)_{H}  \leq ( \lambda^2+1) \nor{v}{H}^2 +\nor{(A - \lambda^2)v }{H}\nor{v}{H} \\
& \leq  ( \lambda^2+2) \nor{v}{H}^2 +\nor{(A - \lambda^2)v }{H}^2.
\end{align*}
Plugging this into~\eqref{e:interm-interm} yields, for all $\mu \geq \mu_0$ and $\lambda \geq 0$,
\begin{align*}
\nor{v}{H}\leq G(\mu) \sqrt{T}\big(  \nor{B^*v}{Y} + C\nor{(A - \lambda^2)v }{H}  \big) +\frac{1}{\mu} \left( \nor{(A - \lambda^2)v }{H}+ (\lambda +\sqrt{2}) \nor{v}{H} \right) .
\end{align*}
We let $\alpha >0$ and choose $\mu = \mu(\lambda) = \max \{\lambda + \sqrt{2} + \alpha, \mu_0 \}$ so that to absorb the last term in the right handside, implying for all $\lambda \geq 0$, 
\begin{align*}
\left(1-\frac{\lambda+\sqrt{2}}{\lambda+\sqrt{2} + \alpha}\right) \nor{v}{H}\leq G(\mu(\lambda)) \sqrt{T} \big(  \nor{B^*v}{Y} + C\nor{(A - \lambda^2)v }{H}  \big) +\frac{1}{\mu(\lambda)} \nor{(A - \lambda^2)v }{H}   .
\end{align*}
We then take $\lambda \geq \mu_0$ so that $\mu(\lambda)= \lambda + \sqrt{2} + \alpha \geq 1$. This implies $\frac{1}{\mu(\lambda)} \nor{(A - \lambda^2)v }{H}  \leq c_0^{-1}G(\mu(\lambda)) \nor{(A - \lambda^2)v }{H}$ and thus, for $\lambda \geq  \mu_0$,
\begin{align*}
\frac{\alpha}{\mu(\lambda)} \nor{v}{H}\leq G(\mu(\lambda)) \sqrt{T} \big(  \nor{B^*v}{Y} + C\nor{(A - \lambda^2)v }{H}  \big) +c_0^{-1}G(\mu(\lambda)) \nor{(A - \lambda^2)v }{H} .
\end{align*}
This concludes the proof of the proposition.
\enp

We finally give a proof of Corollary~\ref{corresolvL}.
\bnp[Proof of Corollary~\ref{corresolvL}]
By assumption, $b\geq \delta>0$ on a non empty open set $\omega$. Since $\M$ is compact, $\sup_{x\in \M} d_\scrL(x,\omega)$ is finite.  
For the hypoelliptic wave equation on $H=Y=L^2(\M)$, we proved in \cite[Theorem~1.15]{LL:17Hypo} that \eqref{th-estimate-k-bis} is satisfied for $A=\scrL$, $B_\omega = B^*_{\omega} =$ multiplication by $\mathds{1}_\omega$, $G(\mu) = Ce^{\kappa \mu^k}$, for all $T> 2\sup_{x\in \M} d_\scrL(x,\omega)$ (where $d_\scrL$ is the appropriate sub-Riemannian distance and $k$ the hypoellipticity index of $\scrL$). Since $\nor{\mathds{1}_{\omega} u}{L^2(\M)}\leq \delta^{-1}\nor{b u}{L^2(\M)}$, the same inequality with different constants remains true with $B=B^*=$ multiplication by $b$. Thus, we deduce from Proposition~\ref{c:res-fct-propres} that~\eqref{e:res-fct-prop-bis} is satisfied (after having fixed $\alpha = 2-\sqrt{2}$) with $\mathsf{G} (\lambda)  =K(1+C) (\lambda +2)  G(\lambda + 2) = C(\lambda +2)  e^{\kappa (\lambda +2)^k}$. 
\enp

\subsection{From free resolvent estimate with an observation term to damped resolvent estimate}
In this section, we start from an estimate for $A$ with an observation term like~\eqref{e:res-fct-prop}, and deduce associated estimates for damped operators.

Now, for later use (see Sections~\ref{s:damped-schro} and~\ref{s:damped-wave} below), we introduce the operators:
\begin{align*}
Q_\lambda & = -i(\A_S - i\lambda ) = A-\lambda  +i BB^* ,\\ 
P_\lambda & = P(i\lambda)= A - \lambda^2  + i \lambda BB^* ,
\end{align*} 
 both with domain $D(Q_\lambda) = D(P_\lambda) = D(A)$.

\begin{proposition}
\label{p:P-lambda}
Let $G_1,G_2 \geq 0$, $\lambda >0$, and $v \in D(A)$, and assume 
\begin{align}
\label{e:res-fct-prop-ter}
\nor{v}{H}\leq G_1 \nor{B^* v}{Y} + G_2 \nor{(A - \lambda^2)v }{H}  .
\end{align}
Then we have
\begin{align}
\label{e:res-fct-prop-ter-P}
\nor{v}{H} & \leq \left(  (G_1\lambda^{-\frac12} + G_2 \sqrt{2} \nor{B}{\L(Y;H)})^2 +2 \sqrt{2} G_2\right)  \nor{P_\lambda v}{H}, \\
\label{e:res-fct-prop-ter-Q}
\nor{v}{H} & \leq \left(  (G_1  + G_2 \sqrt{2} \nor{B}{\L(Y;H)})^2 +2 \sqrt{2} G_2\right)  \nor{Q_{\lambda^2} v}{H}.
\end{align}

In particular, given $\mathsf{G} : \R_+ \to \R_+$ such that $\mathsf{G}(\mu) \geq c_0>0$ on $\R_+$ and $\lambda_0\geq 1$, if~\eqref{e:res-fct-prop-bis} is satisfied, then writing $K = (1+ \sqrt{2} \nor{B}{\L(Y;H)})^2 +  2\sqrt{2} c_0^{-1}$, we have \begin{align}
\label{e:res-fct-prop-terter}
\nor{v}{H} & \leq K \mathsf{G}(|\lambda|)^2 \nor{P_\lambda v}{H} , \quad \text{ for all }v \in D(A)  ,  \lambda \in \R, |\lambda| \geq \lambda_0 , \\
\label{e:res-fct-prop-terter-Q}
\nor{v}{H} & \leq K \mathsf{G}\big(\sqrt{\lambda}\big)^2 \nor{Q_\lambda v}{H} , \quad \text{ for all }v \in D(A)  , \lambda  \geq  \lambda_0^2 .
\end{align}
\end{proposition}
Note that when passing from~\eqref{e:res-fct-prop} to~\eqref{e:res-fct-prop-terter}, we change $\mathsf{G}$ to $\mathsf{G}^2$, which is a loss in general; this is linked to the fact that the proof of Proposition~\ref{p:P-lambda} consists only in a very rough estimate, treating the damping terms $iBB^*$ and $i \lambda BB^*$ as remainders.

\bnp[Proof of Proposition~\ref{p:P-lambda}]
We only prove the result for $P_\lambda$, the analogue proof for $Q_\lambda$ is identical. 

First, we remark that, under the above assumptions, we have
\begin{align}
\label{e:delta-omega}
\lambda \nor{B^*v}{Y}^2  =  \lambda \left( BB^* v , v \right)_{H}  = \Im \left(P_\lambda v , v \right)_H \leq \nor{P_\lambda v}{H}  \nor{v}{H} .
\end{align}
Second, we notice that $(A - \lambda^2)v = P_\lambda v - i\lambda BB^* v$ and thus, using~\eqref{e:delta-omega},
\begin{align*}
 \nor{(A - \lambda^2)v }{H}^2 & \leq  2\nor{P_\lambda v}{H}^2 +2 \lambda \nor{BB^* v}{H}^2 \leq  2\nor{P_\lambda v}{H}^2 +2  \nor{B}{\L(Y;H)}^2 \lambda \nor{B^*v}{Y}^2 \\
 &  \leq  2\nor{P_\lambda v}{H}^2 + 2 \nor{B}{\L(Y;H)}^2 \nor{P_\lambda v}{H}  \nor{v}{H} .
\end{align*}
Plugging the last two estimates in~\eqref{e:res-fct-prop-ter} yields 
\begin{align*}
\nor{v}{H}\leq (G_1 \lambda^{-\frac12} + G_2 \sqrt{2} \nor{B}{\L(Y;H)}) \nor{P_\lambda v}{H}^{\frac12}  \nor{v}{H}^{\frac12}   +G_2 \sqrt{2} \nor{P_\lambda v}{H}  .
\end{align*}
Writing 
\begin{align*}
 (G_1 \lambda^{-\frac12} + G_2\sqrt{2} \nor{B}{\L(Y;H)}) \nor{P_\lambda v}{H}^{\frac12}  \nor{v}{H}^{\frac12} 
 \leq \frac12 (G_1\lambda^{-\frac12} + G_2 \sqrt{2} \nor{B}{\L(Y;H)})^2 \nor{P_\lambda v}{H} +\frac12  \nor{v}{H} ,
\end{align*}
allows to absorb the last term in the left hand-side and implies 
\begin{align*}
\frac12 \nor{v}{H}\leq  \frac12 (G_1\lambda^{-\frac12} + G_2 \sqrt{2} \nor{B}{\L(Y;H)})^2 \nor{P_\lambda v}{H} +G_2 \sqrt{2}  \nor{P_\lambda v}{H}  .
\end{align*}
This concludes the proof of~\eqref{e:res-fct-prop-ter-P}, and~\eqref{e:res-fct-prop-terter} corresponds to the case $G_1=G_2 =\mathsf{G}(\lambda)$. Also, we notice that for $\lambda\in\R$, $\overline{P_{-\lambda}u}=P_{\lambda}\overline{u}$, so the statement in the case $\lambda\geq\lambda_0$ implies the case $\lambda\leq -\lambda_0$. Finally, the proof of~\eqref{e:res-fct-prop-ter-Q} is similar to that of~\eqref{e:res-fct-prop-ter-P} (beware that it should be written for $Q_{\lambda^2}$ and not $Q_{\lambda}$), and~\eqref{e:res-fct-prop-terter-Q} follows from changing $\lambda^2$ into $\lambda$.
\enp
Note that another advantage of Proposition~\ref{p:P-lambda} is that it is flexible enough to support perturbations of the operator $A$ by lower order terms. This was used in \cite{JL:20} where similar estimates were used for application to perturbed operators coming from linearization of a nonlinear equation. See also~\cite{CPSST:19,Burq:19} for recent related perturbation results.

\subsection{Damped Schr\"odinger-type equations}
\label{s:damped-schro}
There are not many references concerning the damped Schr\"odinger equation. So let us start from the beginning. We set $\A_S := i A - BB^*$ with $D(\A_S)=D(A)$, so that~\eqref{eq: damped abstract schro} reformulates as $(\d_t- \A_S )u=0$.

The compact embedding $D(A)\hookrightarrow H$ implies that $\A_S$ has a compact resolvent. First spectral properties of $\A_S$ are described in the following lemma.
\begin{lemma} 
\label{l:spec-AS}
The spectrum of $\A_S$ contains only isolated eigenvalues and we have 
\begin{align}
\nor{(z\id -\A_S)^{-1}}{\L(H)} & \leq \frac{1}{\Re(z)} , \quad \text{ for }\Re(z)>0 , \label{schro-res-1} \\
\nor{(z\id -\A_S)^{-1}}{\L(H)} & \leq \frac{1}{|\Im(z)|} , \quad \text{ for }\Im(z) <0 . \label{schro-res-2}
\end{align}
Moreover, assuming $(Au = zu , B^*u=0 )\implies u=0$, we have 
$$
\Sp(\A_S) \subset   [ - \|B^*\|_{\L(H;Y)}^2 , 0) +  i[0,+\infty) .
$$ 
\end{lemma}
\bnp
The structure of the spectrum comes from the fact that $\A_S$ has a compact resolvent (since so does $A$, and $BB^*$ is bounded). Now, for a general $z \in \C$, we have 
$$
\nor{(z\id -\A_S)u}{H} \nor{u}{H} \geq \Re\left((z\id -\A_S)u , u \right)_{H}  =  \Re(z)\nor{u}{H}^2 + \nor{B^*u}{H}^2  \geq \Re(z)\nor{u}{H}^2 , 
$$
which yields~\eqref{schro-res-1}. The statement~\eqref{schro-res-2} comes from 
$$
\nor{(\A_S-z\id)u}{H} \nor{u}{H} \geq \Im \left((\A_S-z\id)u , u \right)_{H}  =  (A u, u)_H  -\Im(z)\nor{u}{H}^2  \geq -\Im(z)\nor{u}{H}^2 . 
$$
Finally given $z \in \Sp(\A_S)$, there is $u \in D(A)\setminus\{0\}$ such that $\A_S u = zu$. Taking inner product with $u$ yields
$$
z \nor{u}{H}^2 = (\A_S u, u )_{H} = i (A u, u)_H - \nor{B^*u}{H}^2.
$$
In particular, $$
\Re(z) = -\frac{\nor{B^*u}{H}^2}{\nor{u}{H}^2} \in [ - \nor{B^*}{\L(H)}^2 , 0]  , \quad \Im(z) = \frac{(A u, u)_H}{\nor{u}{H}^2} \geq 0 .
$$
Now if $\Re(z) = 0$, this implies $B^* u=0$ and hence $ zu = \A_S u  =i A u$. The assumption then yields $u=0$, which contradicts the fact that $u$ is an eigenvector. Thus $\Sp(\A_S) \cap i\R = \emptyset$.
\enp

We then deduce straightforwardly from Lemmata~\ref{p:P-lambda} and~\ref{l:spec-AS} the following result.
\begin{theorem}
\label{c:Q-lambda-S}
Let  $\mathsf{G} : \R_+ \to \R_+$ be such that $\mathsf{G} (\mu) \geq c_0>0$ on $\R_+$, $\lambda_0\geq 1$, and assume~\eqref{e:res-fct-prop-bis}. Then there exists $K>1$ (the same as in Proposition~\ref{p:P-lambda}), such that
\begin{align*}
& \|(i\lambda \id - \A_S)^{-1}\|_{\L(\H)} \leq K\mathsf{G}\big(\sqrt{\lambda}\big)^2 ,  \quad \text{ for all } \lambda  \geq \lambda_0^2 ,  \\
& \Sp(\A_S) \cap  \Gamma_{\mathsf{G},S}  = \emptyset ,  
\end{align*}
where $\Gamma_{\mathsf{G},S} = \left\{z \in \C , \Im(z)\geq \lambda_0^2 ,  \Re(z) \geq - \frac{1}{K \mathsf{G}\big(\sqrt{\Im(z)}\big)^2 } \right\}$.
Finally, assuming further~\eqref{e:res-fct-prop}, there exists another constant $\widetilde{K}\geq K$ such that
\begin{align*}
& \|(i\lambda \id - \A_S)^{-1}\|_{\L(\H)} \leq \widetilde{K}\mathsf{G}\big(\sqrt{|\lambda|}\big)^2 ,  \quad \text{ for all } \lambda  \in \R, \\
& \Sp(\A_S) \cap  \widetilde{\Gamma}_{\mathsf{G},S}  = \emptyset , 
\end{align*}
where $\widetilde{\Gamma}_{\mathsf{G},S} = \left\{z \in \C ,  \Re(z) \geq - \frac{1}{\widetilde{K} \mathsf{G}\big(\sqrt{|\Im(z)|}\big)^2 } \right\}$.
\end{theorem}
\bnp
The first point is a rewriting of~\eqref{e:res-fct-prop-terter-Q} in Lemma~\ref{p:P-lambda}.

The second point comes from the general fact that $\nor{(z\id-\A_S)^{-1}}{\L(H)}  \geq \frac{1}{\dist(z, \Sp(\A_S))}$ (following from a Neumann series expansion). Hence, we have for $\lambda \geq \lambda_0^2$,
$$
\dist(i\lambda, \Sp(\A_S)) \geq \nor{(i\lambda\id-\A_S)^{-1}}{\L(H)}^{-1} \geq  \left( K  \mathsf{G}\big(\sqrt{\lambda}\big)^2  \right)^{-1} ,
$$
which, together with the localization of the spectrum in Lemma~\ref{l:spec-AS}, proves the second point. 

For the last point, Lemma~\ref{l:spec-AS} ensures that $\lambda\mapsto \|(i\lambda \id - \A_S)^{-1}\|_{\L(\H)} $ is a well defined continuous function on $\R$, which is bounded by $\frac{1}{|\lambda|}$ for $\lambda<0$. On the interval $(-\infty,\lambda_0^2]$, it is therefore bounded by a constant $C_0\leq C_0 c_0^{-2} \mathsf{G}\big(\sqrt{|\lambda|}\big)^2$. This gives the expected estimates for all $\lambda\in \R$ with another $\widetilde{K}=\max \big(K,C_0 c_0^{-2}\big)$. Again, Neumann expansion gives the localization of the spectrum. 
\enp
As a consequence, we deduce the following decay.

\begin{theorem}
\label{t:decay-rates-abstractSchrod}
Let  $\mathsf{G} : \R_+ \to \R_+$ be such that $\mathsf{G} (\mu) \geq c_0>0$ on $\R_+$, $\lambda_0\geq 1$, and assume~\eqref{e:res-fct-prop-bis} and \eqref{e::UCPeigenvalue}. 
Assume further that $\mathsf{G}$ is nondecreasing and set $\mathsf{M}(\lambda) =\mathsf{G}\big(\sqrt{\lambda}\big)^2 $.
Then, there exists $c>0$ such that for all $j \in \N^*$, there is $C_j>0$ such that  for all $u_0 \in D(\A_S^j)$ and associated solution $u$ of \eqref{eq: damped abstract schro}, 
$$
\nor{u(t)}{H} \leq  \frac{C_j}{\mathsf{M}_{\log}^{-1}\left( \frac{t}{c j} \right)^j} \nor{ \A_S^j u_0}{H},  \quad \text{ for all } t>0 ,
$$
where $\mathsf{M}_{\log}$ is defined in~\eqref{Mlog}.
\end{theorem}
Again, $\mathsf{M}_{\log}$ in the result can be replaced by $\mathsf{M}$ if it is polynomial at infinity, according to~\cite[Theorem 2.4]{BT:10}.
\bnp
This is a direct corollary of Theorem~\ref{c:Q-lambda-S} and Theorem~\ref{t:batty-duyckaerts} applied to $\B=\A_S$.
\enp

We may now conclude the proofs of Theorems~\ref{e:decay-schro} and~\ref{t:spectral-schro}.
\bnp[Proof of Theorems~\ref{e:decay-schro} and~\ref{t:spectral-schro}] 

Corollary \ref{corresolvL} implies that \eqref{e:res-fct-prop-bis} is true with $\mathsf{G} (\mu) = C e^{\kappa \mu^k}$. Then, Theorem~\ref{c:Q-lambda-S} implies Theorem~\ref{t:spectral-schro}. Indeed, taking into account~\eqref{e::UCPeigenvalue}, we then obtain that the resolvent is bounded on the positive imaginary axis by a constant times 
$\mathsf{M}(\lambda) =\mathsf{G}\big(\sqrt{\lambda}\big)^2 =  Ce^{2 \kappa^+ \lambda^{k/2}}$ (after having changed the constants slightly).

Finally, we obtain 
$$
\mathsf{M}_{\log}(\lambda) = Ce^{2\kappa^+ \lambda^{k/2}}\left( \log \big(1 + Ce^{2\kappa^+ \lambda^{k/2}}\big) + \log (1+\lambda) \right) \leq Ce^{2\kappa^+ \lambda^{k/2}}
$$
(after having changed the constants slightly), and thus $\mathsf{M}_{\log}^{-1}(t) \geq c \log(t)^{2/k}$ for large $t$. Theorem~\ref{t:decay-rates-abstractSchrod} implies Theorem~\ref{e:decay-schro}.
\enp

\subsection{Damped wave-type equations: semigroup setting and end of the proofs}
\label{s:damped-wave}
Let us now turn this estimate into a resolvent estimate for the generator of the damped wave group, and then into a decay estimate for~\eqref{eq: damped abstract waves}.
We equip $\H=H_1\times H$ with the norm 
$$
\|(u_0,u_1)\|_\H^2 = \|(A+\id)^\frac12 u_0 \|_H^2 + \| u_1 \|_H^2, 
$$
and define the seminorm
$$
|(u_0,u_1)|_\H^2 = \|A^\frac12 u_0 \|_H^2 + \| u_1 \|_H^2.
$$
Of course, if $A$ is coercive on $H$, $|\cdot|_\H$ is a norm on $\H$ equivalent to $\|\cdot \|_\H$.
We define the energy of solutions of \eqref{eq: damped abstract waves} by
$$
E(u(t)) = \frac{1}{2} \big( \|A^\frac12 u\|_H^2 + \|\d_t u\|_H^2 \big)
= \frac12 |(u,\d_t u)|^2_{\H}.
$$
 The damped wave equation~\eqref{eq: damped abstract waves} can be recast on $\H$ as a first order system
\begin{equation}
\label{eq: first order eqation}
\left\{
\begin{array}{l}
\d_t U = \A U , \\
U|_{t=0} = \transp(u_0 , u_1) ,
\end{array}
\right.
\quad 
U = 
\left(
\begin{array}{c}
u \\
\d_t u
\end{array}
\right) , \quad
\A = 
\left(
\begin{array}{cc}
0   &  \id \\
- A & - BB^*
\end{array}
\right)  , \quad
D(\A) =  D(A)\times H_1 .
\end{equation}
The compact embeddings $D(A)\hookrightarrow H_1 \hookrightarrow H$ imply that $D(\A)\hookrightarrow \H$ compactly, and that the operator $\A$ has a compact resolvent. First,  spectral properties of $\A$ are described in the following lemma borrowed from~\cite{Leb:96,AL:14}.
We define the following quadratic family of operator 
\begin{align}
\label{e:def-Pz}
P(z) = A + z^2 \id + z BB^* , \quad z \in \C, \quad  D(P(z)) = D(A) .
\end{align}
\begin{lemma}[Lemma~4.2 of~\cite{AL:14}]
\label{lemma: check assumptions}
The spectrum of $\A$ contains only isolated eigenvalues and, provided \eqref{e::UCPeigenvalue} is satisfied, we have 
$$
\Sp(\A) \subset \left( \big( - \frac12 \|B^*\|_{\L(H;Y)}^2 , 0 \big) + i\R\right) \cup 
\left([ - \|B^*\|_{\L(H;Y)}^2 , 0] + 0 i \right)  ,
$$
with $\ker(\A) = \ker(A)\times \{0\}$.
Moreover, the operator $P(z)$ in~\eqref{e:def-Pz} is an isomorphism from $D(A)$ onto $H$ if and only if $z \notin \Sp(\A)$.
\end{lemma}
This Lemma leads us to introduce the spectral projector of $\A$ on $\ker(\A)$, given by 
$$
\Pi_0 = \frac{1}{2i\pi}\int_{\gamma} (z\id - \A)^{-1} dz \in \L(\H),
$$
where $\gamma$ denotes a positively oriented circle centered on $0$ with a radius so small that $0$ is the single eigenvalue of $\A$ in the interior of $\gamma$.
We set $\dot{\H} = (\id - \Pi_0)\H$ and equip this space with the norm
$$
\|(u_0 , u_1)\|_{\dot{\H}}^2 := |(u_0 , u_1)|_{\H}^2 
= \|A^{\frac12}u_0\|_{H}^2 + \|u_1 \|_{H}^2 ,
$$
and associated inner product.
This is indeed a norm on $\dot{\H}$ since $\|(u_0 , u_1)\|_{\dot{\H}} =0$ is equivalent to $(u_0 , u_1)\in \ker(A)\times \{0\} = \Pi_0\H$.
Besides, we set $\dot{\A} = \A|_{\dot{\H}}$ with domain $D(\dot{\A}) = D(\A) \cap \dot{\H}$. 
Remark that $\Sp(\dot{\A}) = \Sp(\A) \setminus \{0\}$ and thus $\Sp(\dot{\A}) \cap i\R = \emptyset$.

\begin{lemma}[Lemma~4.3 of~\cite{AL:14}]
\label{lemma: semigroups}
The operator $\dot{\A}$ generates a contraction $C^0$-semigroup on~$\dot{\H}$, denoted $(e^{t\dot{\A}})_{t\geq 0}$. Moreover, the operator $\A$ generates a bounded $C^0$-semigroup on~$\H$, denoted $(e^{t\A})_{t\geq 0}$ and the unique solution to~\eqref{eq: damped abstract waves} is given by $(u, \d_t u )(t) = e^{t\A}(u_0,u_1)$. Finally,  we have
\begin{equation}
\label{eq: decomposition semigroup}
e^{t\A}= e^{t\dot{\A}} (\id -\Pi_0)  + \Pi_0  , \quad \text{for all } t\geq 0 .
\end{equation}
\end{lemma} 

Once we have put the abstract damped wave equation~\eqref{eq: damped abstract waves} in the appropriate semigroup setting, it remains to: 
\begin{enumerate}
\item \label{point-1} deduce from~\eqref{e:res-fct-prop}-\eqref{e:res-fct-prop-bis} a resolvent estimate for $\dot{\A}$,
\item \label{point-2} relate this resolvent estimate to a decay estimate for $e^{t\dot{\A}}$, and
\item \label{point-3} deduce the decay of the energy for~\eqref{eq: damped abstract waves}.
\end{enumerate}

\medskip
Step~\ref{point-1} is realized thanks to the following result from~\cite{AL:14}.
\begin{lemma}[Lemma~4.6 of~\cite{AL:14}]
\label{lemma: conditions stable resolvent}
There exist $C>1$ such that for $s \in \R$, $|s|\geq 1$, 
\begin{align}
\label{eq: estimate A estimate dotA}
 C^{-1}\|(is\id - \dot{\A})^{-1}\|_{\L(\dot{\H})} - \frac{C}{|s|} 
\leq \|(is\id - \A)^{-1}\|_{\L(\H)} 
\leq  C\|(is\id - \dot{\A})^{-1}\|_{\L(\dot{\H})} + \frac{C}{|s|}  ,
\end{align}
\begin{align}
\label{eq: estimate P equiv A}
C^{-1}|s| \|P(is)^{-1}\|_{\L(H)} \leq  \|(is\id - \A)^{-1}\|_{\L(\H)}
\leq C \left( 1 + |s|\|P(is)^{-1}\|_{\L(H)} \right).
\end{align}
\end{lemma}
As a corollary of this together with Proposition~\ref{p:P-lambda}, we deduce the following result.

\begin{theorem}
\label{c:P-lambda-A}
Let  $\mathsf{G} : \R_+ \to \R_+$ be such that $\mathsf{G} (\mu) \geq c_0>0$ on $\R_+$, $\lambda_0\geq 1$, and assume~\eqref{e:res-fct-prop-bis}. Then there exists $K>1$ such that
\begin{align*}
& \|(i\lambda \id - \A)^{-1}\|_{\L(\H)} \leq K|\lambda| \mathsf{G}(|\lambda|)^2 ,  \quad \text{ for all } \lambda\in \R, |\lambda|  \geq \lambda_0 ,  \\
& \|(is\id - \dot{\A})^{-1}\|_{\L(\dot{\H})}\leq K|\lambda|  \mathsf{G}(|\lambda|)^2  ,  \quad \text{ for all } \lambda\in \R, |\lambda|  \geq \lambda_0  , \\
&\Sp(\dot{\A}) \cap \Gamma_{\mathsf{G}} =\emptyset , \qquad  \Sp(\A) \cap  \Gamma_{\mathsf{G}}  = \emptyset , 
\end{align*}
where $\Gamma_{\mathsf{G}} = \left\{z \in \C , |\Im(z)| \geq \lambda_0 ,  \Re(z) \geq - \frac{1}{K |\Im(z)|  \mathsf{G}(|\Im(z)|)^2} \right\}$.

Finally, assuming further~\eqref{e:res-fct-prop}, there exists another constant $\widetilde{K}\geq K$ such that
\begin{align*}
& \|(is\id - \dot{\A})^{-1}\|_{\L(\dot{\H})}\leq \widetilde{K}\left\langle \lambda\right\rangle  \mathsf{G}(|\lambda|)^2  ,  \quad \text{ for all } \lambda\in \R, \\
&\Sp(\dot{\A}) \cap \widetilde{\Gamma}_{\mathsf{G}} =\emptyset , \qquad  \Sp(\A) \cap  \Gamma_{\mathsf{G}}  = \{0\}, 
\end{align*}
where $\widetilde{\Gamma}_{\mathsf{G}} = \left\{z \in \C , \Re(z) \geq - \frac{1}{\widetilde{K}\left\langle \Im(z) \right\rangle \mathsf{G}( |\Im(z)|)^2} \right\}$.
\end{theorem}

\bnp[Proof of Theorem~\ref{c:P-lambda-A}]
The first two points are corollaries of~\eqref{e:res-fct-prop-terter} in Proposition~\ref{p:P-lambda} combined with Lemma~\ref{lemma: conditions stable resolvent}.
 
The last point comes from $\Sp(\dot{\A}) = \Sp(\A) \setminus \{0\}$, together with the general fact that $\nor{(z\id-\dot{\A})^{-1}}{\L(\H)}  \geq \frac{1}{\dist(z, \Sp(\dot{\A}))}$ (following from a Neumann series expansion). Hence, we have for $\lambda\in \R$, $|\lambda| \geq \lambda_0$,
$$
\dist(i\lambda, \Sp(\dot{\A})) \geq \nor{(i\lambda\id-\dot{\A})^{-1}}{\L(\H)}^{-1} \geq  \left( K |\lambda|  \mathsf{G}(|\lambda|)^2  \right)^{-1} ,
$$
which, together with the localization of the spectrum in Lemma~\ref{lemma: check assumptions}, proves the statement about the zone free of spectrum. The proof concerning the compact zone follows the same way as Theorem \ref{c:Q-lambda-S}.
\enp

\medskip
Step~\ref{point-2} is achieved as a consequence of Theorem~\ref{t:batty-duyckaerts} applied to the operator $\B=\dot{\A}$.

\medskip
Finally, Step~\ref{point-3} is a consequence of the following elementary lemma~\ref{l:AL:14-44}, linking the energy of solutions to the abstract damped wave equation~\eqref{eq: damped abstract waves} to the norm of the semigroup $\big(e^{t\dot{\A}}\big)_{t\geq0}$.
\begin{lemma}
\label{l:AL:14-44}
For all $j \in \N^*$, $U_0 \in D(\A^j)$ such that $\Pi_0 U_0 \neq U_0$, and associated solution $u$ of \eqref{eq: damped abstract waves}, we have
\begin{align*}
\frac{E(u(t))}{\frac12 | \A^j U_0 |_{\H}^2} = \frac{|e^{t\A} U_0 |_{\H}^2}{| \A^j U_0 |_{\H}^2}  =\frac{\|e^{t\dot{\A}}\dot{U}_0\|_{\dot{\H}}^2}{\|\dot{\A}^j \dot{U}_0\|_{\dot{\H}}^2} , \quad \text{where} \quad  \dot{U}_0 = (\id-\Pi_0)U_0 .
\end{align*}
In particular, setting $f_j(t) := \nor{e^{t\dot{\A}}\dot{\A}^{-j}}{\L(\dot{\H})}$ for $j \in \N^*$, we have for all $U_0 \in D(\A^j)$ and associated solution $u$ of \eqref{eq: damped abstract waves},
$$
E(u(t)) \leq \frac12 f_j(t)^2\| \A^j U_0 \|_{\H}^2 , \quad \text{ for all }t\geq0 .
$$
\end{lemma}
\bnp
This is essentially~\cite[Lemma~4.4]{AL:14}. Recalling that $\A U_0= \dot{\A} \dot{U}_0$, we have
\begin{align*}
E(u(t)) & = \frac{1}{2} \big( \|A^\frac12 u(t)\|_H^2 + \|\d_t u(t)\|_H^2 \big) =\frac12 |e^{t\A} U_0 |_{\H}^2 
= \frac12 |e^{t\dot{\A}} \dot{U}_0 + \Pi_0 U_0|_{\H}^2
= \frac12 \|e^{t\dot{\A}}\dot{U}_0\|_{\dot{\H}}^2 ,  \\
\|\dot{\A}^j\dot{U}_0\|_{\dot{\H}}^2 & = | \A^j U_0 |_{\H}^2 =  \| \A^j U_0 \|_{\H}^2 ,
\end{align*}
which yields the first statement. The second one follows from the fact that $|\cdot|_{\H} \leq \nor{\cdot}{\H}$.
\enp

As a consequence, we deduce the following decay.

\begin{theorem}
\label{t:decay-rates-abstract}
Let  $\mathsf{G} : \R_+ \to \R_+$ be such that $\mathsf{G} (\mu) \geq c_0>0$ on $\R_+$, $\lambda_0\geq 1$, and assume~\eqref{e::UCPeigenvalue} and \eqref{e:res-fct-prop-bis}.
Assume further that $\mathsf{G}$ is nondecreasing and set $\mathsf{M}(\lambda) = \left\langle \lambda\right\rangle \mathsf{G}(\lambda)^2$.
Then, there exists $c>0$ such that for all $j \in \N^*$, there is $C_j>0$ such that  for all $U_0 \in D(\A^j)$ and associated solution $u$ of \eqref{eq: damped abstract waves}, 
$$
E(u(t))^\frac12\leq  \frac{C_j}{\mathsf{M}_{\log}^{-1}\left( \frac{t}{c j} \right)^j} \nor{ \A^j U_0}{\H},  \quad \text{ for all } t>0 ,
$$
where $\mathsf{M}_{\log}$ is defined in~\eqref{Mlog}.
\end{theorem}
Again, $\mathsf{M}_{\log}$ in the result can be replaced by $\mathsf{M}$ if it is polynomial at infinity, according to~\cite[Theorem 2.4]{BT:10}.
\bnp
This is a direct corollary of Theorem~\ref{c:P-lambda-A}, Theorem~\ref{t:batty-duyckaerts} applied to $\mathcal{X} = \dot{\H}$ and $\B=\dot{\A}$, together with Lemma~\ref{l:AL:14-44}.
\enp
We conclude this paragraph with the proofs of Theorems \ref{e:decay} and \ref{t:spectral}.
\bnp[Proof of Theorems \ref{e:decay} and \ref{t:spectral}]
Again, Corollary \ref{corresolvL} implies the unique continuation property \eqref{e:UCP-eig-Hypo} (that is \eqref{e::UCPeigenvalue} in the present context) together with~\eqref{e:res-fct-prop-bis} with $\mathsf{G} (\mu) = C e^{\kappa \mu^k}$. With this estimate at hand, Theorem~\ref{e:decay} is an application of Theorem \ref{t:decay-rates-abstract} with $\mathsf{M}(\lambda) = \left\langle \lambda\right\rangle \mathsf{G}(\lambda)^2\leq  Ce^{2 \kappa^+ \lambda^{k}}$ (after having changed the constants slightly), while Theorem \ref{t:spectral} is implied by Lemma \ref{lemma: check assumptions} and Theorem \ref{c:P-lambda-A}. 
\enp

\subsection{Damped plate equation}
The plate equation actually fits into the ``wave-type'' framework. Indeed, the abstract plate equation
\begin{align}
\label{eq: damped abstract beam}
\begin{cases}
\d_t^2 u + A^2 u + B B^* \d_t u = 0, \\
(u, \d_t u)|_{t=0} = (u_0 , u_1) , 
\end{cases}
\end{align} 
is actually a particular case of the abstract equation \eqref{eq: damped abstract waves} applied with the operator $A^2$ (instead of $A$) which is still nonnegative selfadjoint with compact resolvent. In this case, we define $H_2 = D(A)$, equipped with the graph norm $\nor{u}{H_2} := \|(A^2+\id)^\frac12 u \|_H$, and its dual $H_{-2} = (H_2)'$ (using $H$ as a pivot space) endowed with the norm $\nor{u}{H_{-2}}:=\|(A^2+\id)^{-\frac12} u \|_H$.

The natural space is then $\H=H_2\times H$ with the norm 
$$
\|(u_0,u_1)\|_\H^2 = \|(A^2+\id)^\frac12 u_0 \|_H^2 + \| u_1 \|_H^2, 
$$
and the seminorm
$$
|(u_0,u_1)|_\H^2 = \|A u_0 \|_H^2 + \| u_1 \|_H^2.
$$
The associated energy is 
$$
E_P(u(t)) = \frac{1}{2} \big( \|A u\|_H^2 + \|\d_t u\|_H^2 \big)
= \frac12 |(u,\d_t u)|^2_{\H}.
$$
In order to transfer the properties of $A$ to $A^2$, we will only need the following simple lemma.
 \begin{lemma}
Assume \eqref{e:res-fct-prop-bis} is satisfied. Then, we have
\begin{align}
\label{e:res-fct-prop-beam}
\nor{v}{H}\leq \mathsf{G}(\sqrt{\lambda}) \big( \nor{B^* v}{Y} +\lambda^{-1} \nor{(A^2 - \lambda^2)v }{H}  \big) , \quad \text{ for all }v \in D(A^2) , \lambda \geq \lambda_0^2 .
\end{align}
\end{lemma}
\bnp
Since $A$ is a nonnegative operator, we have $\nor{(A + \lambda^2)w }{H}\geq  \lambda^2 \nor{w }{H}$ for all $w\in D(A)$. Applying this to $w=(A - \lambda^2)v$ gives $\nor{(A^2 - \lambda^4)v }{H}\geq \lambda^2 \nor{(A - \lambda^2)v}{H}$.
This, combined with \eqref{e:res-fct-prop-bis} implies
\begin{align}
\nor{v}{H}\leq \mathsf{G}(\lambda) \big( \nor{B^* v}{Y} + \nor{(A - \lambda^2)v }{H}  \big)\leq \mathsf{G}(\lambda) \left( \nor{B^* v}{Y} + \frac{1}{\lambda^2}\nor{(A^2 - \lambda^4)v }{H}  \right) .
\end{align}
This is the expected result up to changing $\lambda$ to $\sqrt{\lambda}$. 
\enp
The previous Lemma implies that if \eqref{e:res-fct-prop-bis} is satisfied, the assumptions of Theorem \ref{t:decay-rates-abstract} are satisfied for the operator $A^2$ with $G_P(\lambda)=G(\sqrt{\lambda})$. It directly gives the following result.
\begin{theorem}
\label{t:decay-rates-abstractbeam}
Let  $\mathsf{G} : \R_+ \to \R_+$ be such that $\mathsf{G} (\mu) \geq c_0>0$ on $\R_+$, $\lambda_0\geq 1$, and assume~\eqref{e::UCPeigenvalue} and~\eqref{e:res-fct-prop-bis}. 
Assume further that $\mathsf{G}$ is nondecreasing and set $\mathsf{M}(\lambda) = \lambda \mathsf{G}\big(\sqrt{\lambda}\big)^2$.
Then, there exists $c>0$ such that for all $j \in \N^*$, there is $C_j>0$ such that  for all $U_0 \in D(\A^j)$ and associated solution $u$ of~\eqref{eq: damped abstract beam},
$$
E_P(u(t))^\frac12\leq  \frac{C_j}{\mathsf{M}_{\log}^{-1}\left( \frac{t}{c j} \right)^j} \nor{ \A_P^j U_0}{\H},  \quad \text{ for all } t>0 ,
$$
where $\mathsf{M}_{\log}$ is defined in~\eqref{Mlog}.
\end{theorem}

\bnp[Proof of Theorem \ref{e:decaybeam}]
Thanks to Corollary \ref{corresolvL}, \eqref{e:res-fct-prop-bis} is true with $\mathsf{G} (\mu) = C(\mu +2)  e^{\kappa (\mu +2)^k}$. Theorem \ref{e:decaybeam} is an application of Theorem \ref{t:decay-rates-abstractbeam} with $\mathsf{M}(\lambda) = \lambda \mathsf{G}(\sqrt{\lambda})^2\leq  Ce^{2 \kappa^+ \lambda^{k/2}}$ (after having changed the constants slightly).
\enp

 \subsection{Lower bounds: proof of Proposition~\ref{Prop:BCG-damped}}

\bnp[Proof of Proposition~\ref{Prop:BCG-damped}]
According to~\cite[Proposition~1.14]{LL:17Hypo} (which relies on~\cite[Section~2.3]{BeauchardCanGugl}), since $\supp(b) \cap \left\{x_1=0\right\} = \emptyset$, there exist $C,c_0>0$ and a sequence $(\lambda_j , \varphi_j) \in \R_+ \times C^\infty(\M)$ such that 
$$
\scrL\varphi_j = \lambda_j \varphi_j ,\quad  \varphi_j|_{\d\M} = 0 , \quad \|\varphi_j\|_{L^2(\M)}=1, \quad \lambda_j \to +\infty , \quad \| \varphi_j\|_{L^2(\supp(b))} \leq  Ce^{ - c_0 \lambda_j^{\frac{k}{2}}} .
$$
As a consequence, concerning the damped Schr\"odinger resolvent, we have 
$$
\nor{(\A_S - i\lambda_j) \varphi_j}{L^2(\M)} = \nor{(i\scrL - b - i\lambda_j) \varphi_j}{L^2(\M)} = \nor{ b \varphi_j}{L^2(\M)}  \leq \nor{b}{L^\infty}Ce^{ - c_0 \lambda_j^{\frac{k}{2}}}  .
$$
This implies the second estimate in~\eqref{e:lower-resolvent-bounds} with $s_j = \lambda_j$.

Concerning the damped wave resolvent, recalling the definition of $P(z)$ in~\eqref{e:def-Pz}, we write 
$$
\nor{P\big(i\sqrt{\lambda_j}\big)\varphi_j}{L^2} = \nor{\big(\scrL - \lambda_j + i\sqrt{\lambda_j} b\big)\varphi_j}{L^2} = \nor{\sqrt{\lambda_j} b\varphi_j}{L^2} \leq  \sqrt{\lambda_j} \nor{b}{L^\infty}Ce^{ - c_0 \lambda_j^{\frac{k}{2}}} .
$$
With $s_j = \sqrt{\lambda_j}$, this implies $\nor{P\big(is_j\big)\varphi_j}{L^2} \leq s_j Ce^{ - c_0 s_j^{k}}$, and using~\eqref{eq: estimate P equiv A} in Lemma~\ref{lemma: conditions stable resolvent} proves the first estimate in~\eqref{e:lower-resolvent-bounds}.

\medskip
The last part of the Proposition follows from~\eqref{e:lower-resolvent-bounds} together with the first implication in Theorem~\ref{t:batty-duyckaerts} (and, in case of damped waves, equivalence between the resolvents of $\A$ et $\dot{\A}$ in~\eqref{eq: estimate A estimate dotA} in Lemma~\ref{lemma: conditions stable resolvent}).
\enp

\small
\bibliographystyle{alpha}
\bibliography{bibli}
\end{document}

%% file: mymacros.tex
\newtheorem{lemma}{Lemma}[section]
\newtheorem{theorem}[lemma]{Theorem}
\newtheorem{proposition}[lemma]{Proposition}
\newtheorem{corollary}[lemma]{Corollary}

\theoremstyle{definition}

\newtheorem{definition}[lemma]{Definition}

\makeatletter
\def\keywords{
    \vspace{1ex}
    \noindent
    \if@twocolumn
      \small{\bf  Keywords}\/---$\!$    \else
      \begin{center}\small\ {\bf Keywords}\end{center}\quotation\small
    \fi}
\def\endkeywords{\vspace{0.6em}\par\if@twocolumn\else\endquotation\fi
    \normalsize\rm}
\makeatother

\renewcommand{\O}{\ensuremath{\mathcal O}}

\renewcommand{\L}{\ensuremath{\mathcal L}}

\newcommand{\scrL}{\ensuremath{\mathscr L}}

\DeclareMathOperator{\Sp}{Sp}
\DeclareMathOperator{\loc}{loc}

\DeclareMathOperator{\Lie}{Lie}

\newcommand{\mb}[1]{\ensuremath{\mathbb{#1}}}
\newcommand{\N}{{\mb{N}}}

\newcommand{\R}{{\mb{R}}}
\newcommand{\C}{{\mb{C}}}

\newcommand{\eps}{\varepsilon}

\newcommand{\M}{\ensuremath{\mathcal M}}
\newcommand{\D}{\ensuremath{\mathscr D}}
\newcommand{\B}{\ensuremath{\mathcal B}}
\newcommand{\A}{\ensuremath{\mathcal A}}

\let \Re \relax
\DeclareMathOperator{\Re}{Re}
\let \Im \relax
\DeclareMathOperator{\Im}{Im}

\newcommand{\ovl}[1]{\overline{#1}}

\newcommand{\Con}{\ensuremath{\mathscr C}}

\DeclareMathOperator{\supp}{supp}

\DeclareMathOperator{\dist}{dist}

\DeclareMathOperator{\vect}{span}

\DeclareMathOperator{\id}{Id}

\newcommand{\transp}{\ensuremath{\phantom{}^{t}}}

\renewcommand{\d}{\ensuremath{\partial}}

\newcommand{\Z}{\mathbb Z}

\renewcommand{\H}{\ensuremath{\mathcal H}}

